\documentclass[11pt]{article}
\usepackage{amssymb} 
\newtheorem{lem}{Lemma} 
\newtheorem{thm}{Theorem}

\begin{document}

\title{Chords of 2-factors in planar cubic bridgeless graphs}

\author{Ajit A. Diwan \\
Department of Computer Science and Engineering,\\
Indian Institute of Technology Bombay, Mumbai 400076, India.\\
email:\texttt{aad@cse.iitb.ac.in}}

\maketitle

\begin{abstract}
We show that every edge in a 2-edge-connected planar cubic graph is either contained in a 2-edge-cut
or is a chord of some cycle that is contained in a 2-factor of the graph. As a consequence, we show
that every edge in a cyclically 4-edge-connected planar cubic graph, except $K_2^3$ and $K_4$, is contained 
in a perfect matching whose removal disconnects the graph. We obtain a complete characterization of 
2-edge-connected planar cubic graphs that have an edge such that every 2-factor containing the edge is a 
Hamiltonian cycle, and also of those that have an edge such that the complement of every perfect matching 
containing the edge is a Hamiltonian cycle. Another immediate consequence of the main result is that for any 
two edges contained in a facial cycle of a 2-edge-connected planar cubic graph, there exists a 2-factor
in the graph such that both edges are contained in the same cycle of the 2-factor. We conjecture
that this property holds for any two edges in a 2-edge-connected planar cubic graph, and prove it
for planar cubic bipartite graphs. The main result is proved in the dual form by showing that every plane 
triangulation admits a vertex 3-coloring such that no face is monochromatic and there is exactly one
specified edge between a specified pair of color classes.
\end{abstract}

\section{Introduction}
A classical result in graph theory is Petersen's theorem~\cite{P} that every 2-edge-connected 
cubic graph has a perfect matching, and hence a 2-factor obtained by taking the complement of the matching. 
This can be strengthened to show that there exists a perfect matching including any specified edge, and a 
2-factor including any two specified edges in the graph~\cite{S}. More can be said about the structure of the 
2-factor if the cubic graph has additional properties, in particular, for 2-edge-connected planar cubic graphs. 
It is well-known that the four color theorem is equivalent to the statement that every 2-edge-connected planar 
cubic graph is 3-edge-colorable. This is equivalent to the statement that every such graph has a 2-factor in 
which all cycles have even lengths. It was shown in~\cite{D} that every 2-edge-connected planar cubic graph
with at least 6 vertices has a \emph{disconnected} 2-factor, that is, a 2-factor with more than one
component. H{\"a}ggkvist~\cite{H} showed that every edge in a \emph{bipartite} cubic graph is either contained
in a 2-edge-cut or is a chord of some cycle that is contained in a 2-factor of the graph. Barnette~\cite{B}
conjectured that every 3-edge-connected planar cubic bipartite graph has a Hamiltonian cycle, which is a 
connected 2-factor, and this question is still open.

In this paper, we show that H{\"a}ggkvist's result also holds for all 2-edge-connected planar cubic graphs.
In other words, every edge in a 2-edge-connected planar cubic graph is either contained in a 2-edge-cut
or is a chord of some cycle that is contained in a 2-factor in the graph.  An edge in a cubic graph is a chord 
of some cycle iff it is not contained in an edge-cut of size at most two. Thus another way of saying this is 
that if an edge is a chord of some cycle then it is a chord of some cycle that is contained in a 2-factor. As a 
consequence, we get a simpler and completely different proof of the result in~\cite{D} and in fact prove much 
stronger statements. We show that every edge in a cyclically 4-edge-connected planar cubic graph, except 
$K_2^3$ and $K_4$, is contained in a perfect matching whose complement is a disconnected 2-factor. We call such 
a perfect matching a \emph{separating} perfect matching. We obtain a complete characterization of 
2-edge-connected planar cubic graphs that have an edge that is not contained in a separating perfect matching, 
and also of those 2-edge-connected planar cubic graphs that have an edge such that every 2-factor containing 
the edge is a Hamiltonian cycle. It may be noted that a characterization of cubic bipartite graphs in which 
every 2-factor is a Hamiltonian cycle has been conjectured~\cite{F}, but the problem is still open.

Another consequence of our result is that for any two edges contained in a facial cycle of a 
2-edge-connected planar cubic graph, there exists a cycle containing both the edges that is
contained in a 2-factor of the graph. We conjecture that this property holds for any
two edges in the graph. It is also possible that it holds for any two edges in a cubic bipartite
graph. We prove it for any two edges in a planar cubic bipartite graph.

The Petersen graph, which is the smallest 2-edge-connected cubic graph that is not 3-edge-colorable,
also shows that in general, a 2-edge-connected cubic graph may not have any 2-factor in which
some cycle has a chord. In this case, no edge is contained in a 2-edge-cut and no edge is a chord of
a cycle contained in a 2-factor. It is an interesting question to determine which 2-edge-connected cubic 
graphs have a 2-factor in which some cycle has a chord. It is possible that every 3-edge-colorable cubic 
graph has a 2-factor in which some cycle has a chord, although it is not true that every edge that is not 
contained in a 2-edge-cut is a chord of some cycle contained in a 2-factor in such a graph. It has been 
conjectured~\cite{OP} that every 3-edge-colorable cubic graph has an edge $e$ such that the graph obtained by 
deleting $e$ and contracting one other edge incident with each endpoint of $e$ is also 3-edge-colorable. 
Perhaps a stronger statement that combines both properties may hold. If every 3-edge-colorable cubic graph has 
a 3-edge-coloring such that some 2-edge-colored cycle has a chord, it would imply both the results. This 
statement can be verified easily for planar cubic graphs, even without assuming that all of them are
3-edge-colorable. It also holds for cubic bipartite graphs by H{\"a}ggkvist's result, since any such graph 
must have an edge not contained in a 2-edge-cut.

One reason for studying 2-factors in 2-edge-connected planar cubic graphs is that results on these have
a dual version in terms of vertex colorings of the dual graphs. The planar duals of 2-edge-connected
planar cubic graphs are plane triangulations, which are loopless plane graphs embedded in the plane so
that each face has 3 edges on its boundary. The most famous example of this is of course Tait's
reformulation of the four color theorem. Tait in fact falsely conjectured that all 3-edge-connected
planar cubic graphs are Hamiltonian, which would have proved the four color theorem. Penaud~\cite{Pe} 
first noted that the existence of a 2-factor in a 2-edge-connected planar cubic graph implies that any 
plane triangulation has a non-monochromatic 2-coloring, that is a coloring of the vertices with 2 
colors such that no three vertices that are on the boundary of the same face have the same color.  As shown 
in~\cite{DK}, the existence of a disconnected 2-factor in a 2-edge-connected planar cubic graph  implies the 
existence of a strict 3-coloring of the vertices of a plane triangulation such that no face is monochromatic 
and also no face is rainbow, that is, some two vertices in the boundary of any face have the same color.
In general, it is shown in~\cite{DK} that the existence of a 2-factor with at least $k$ components
is equivalent to the existence of a strict $(k+1)$-coloring of the dual such that no face is monochromatic 
or rainbow. Although this problem is NP-complete for arbitrary $k$, it is still open for any fixed $k \ge 3$.

The dual version of finding a 2-factor with a chord is to find a non-monochromatic
2-coloring of the vertices of a triangulation such that the subgraph formed by the monochromatic edges
has a bridge. This is equivalent to finding a non-monochromatic 3-coloring such that for some two 
color classes there is exactly one edge (the bridge) between the color classes. To find a 2-factor
in which a specified edge is a chord of some cycle, we fix the edge between the color classes to be
the dual of the edge required to be a chord. We will use this dual formulation and prove the existence
of such a 3-coloring by induction. The dual result may also be of independent interest.

\section{Terminology}

The notation and terminology used is mostly standard. We will only clarify the terms that are specific
to this work. We consider undirected graphs that may have multiple edges but no self-loops. A graph is 
$k$-edge-connected if it cannot be disconnected by removing less than $k$ edges. An edge whose removal 
increases the number of connected components in a graph is called a bridge. A graph is cubic if the degree 
of every vertex is three. The cubic graph with two vertices and three edges joining them is denoted by $K_2^3$.
A cubic graph is said to be cyclically $k$-edge-connected if removing any set of less than $k$ edges
results in a graph with at most one component containing a cycle. A plane graph is a graph that has been 
embedded in the plane. Any such embedding divides the plane into connected regions called faces. The unbounded 
region will be called the external face and all other faces are said to be internal. We will identify a face by 
the sequence of vertices and edges on its boundary in circular order. A plane triangulation is a plane graph 
such that every face is bounded by 3 edges. The dual of any plane cubic graph is a plane triangulation and
vice-versa. Since we consider only bridgeless cubic graphs, their duals will be triangulations without any 
self-loops. This implies the plane triangulations are necessarily 2-connected and no vertex or edge is repeated 
on the boundary of any face.  A plane near-triangulation is a 2-connected plane graph such that every
internal face is bounded by 3 edges. The boundary of a near-triangulation is the sequence of vertices that 
are on the boundary of the external face, in circular order. The vertices and edges that are on the
boundary of the external face are called boundary vertices and edges, respectively, and the other vertices 
and edges are said to be internal. A chord of a near-triangulation is an internal edge both of whose endpoints 
are boundary vertices.

A $3$-coloring of a graph $G$ is a function $f : V(G) \rightarrow \{a,b,c\}$ that assigns
one of 3 colors to each vertex in $G$.  An internal face of a near-triangulation is said to be
monochromatic in a 3-coloring $f$ if all 3 vertices on the boundary of the face have the same color.
The coloring $f$ is said to be non-monochromatic if there is no monochromatic internal face in $f$.
The coloring is said to be a 2-coloring with colors $\{a,c\}$ if $f(v) \neq b$ for all $v \in V(G)$.
 
A non-monochromatic 3-coloring $f$ of a near-triangulation $G$  with boundary $v_1,v_2,\ldots,v_l$ is said 
to be \emph{special} if for any two vertices $u,v$ such that $f(u) = a$ and $f(v) = b$, $u$ is adjacent to $v$ 
iff $u = v_1$ and $v = v_2$. An adjacent pair of vertices $u,v$ is called an $ab$-pair in a coloring $f$ if 
$f(u) = a$ and $f(v) = b$.  Thus a 3-coloring of $G$ is special if it is non-monochromatic and the only 
$ab$-pair is $v_1,v_2$. For a string $s = s_1s_2\ldots s_l$ of length $l$ over the alphabet 
$\{a,b,c\}$, we say that $G$ has a special 3-coloring $f$ with colors $s$ assigned to the boundary if $f$ is a 
special 3-coloring of $G$ such that $f(v_i) = s_i$ for $1 \le i \le l$. We say that $s$ is \emph{feasible} for 
$G$ if there exists a special 3-coloring of $G$ with colors $s$ assigned to the boundary. Note that we can 
choose the vertex $v_1$ in the boundary arbitrarily and label the other vertices in circular order, either  
clockwise or anti-clockwise. In the case of a triangulation, we can choose the unbounded face also 
arbitrarily, so that any pair of adjacent vertices in a triangulation may be labeled $v_1,v_2$. 

We now state some basic results using the above terminology.

\begin{lem}[Penaud]
\label{Pet}  
Let $G$ be a plane triangulation without self-loops having the boundary $v_1,v_2,v_3$. Then there exists
a non-monochromatic 2-coloring of $G$ with colors $\{a,c\}$ such that the boundary is assigned
colors $aac$.
\end{lem}

\noindent {\bf Proof:}
This follows from the stronger form of Petersen's theorem. Consider a matching in the dual that
contains the dual edge of $v_1v_2$. Deleting the edges in this matching gives a 2-factor which is a 
collection of disjoint cycles in the plane. The regions into which the plane is divided by these cycles
can be 2 colored so that adjacent regions get distinct colors. Assign a color to a region depending on the
parity of the number of cycles in the 2-factor whose interior contains the region. Assign a vertex in the 
triangulation the color of the region containing the face in the dual corresponding to the vertex. Then no 
face in the triangulation will be monochromatic since the vertex in the dual corresponding to the face will
be contained in some cycle in the 2-factor, and all 3 faces incident with the vertex cannot be on the
same side of the cycle. However, they will be on the same side of any other cycle, which implies all three
cannot have the same color. Since the dual edge of $v_1v_2$ is in the matching, $v_1$ and $v_2$
will get the same color. We can assume this color is $a$, without loss of generality, and since $v_1,v_2,v_3$
is a face in $G$, $v_3$ must have a different color.
\hfill $\Box$

A small modification of this gives the coloring equivalent of a 2-factor with a chord.

\begin{lem}
\label{3col}
Let $G$ be a 2-edge-connected plane cubic graph and $uv$ an edge in $G$ that is not contained in
a 2-edge-cut. Let $G'$ be the dual of $G$ with the endpoints of the dual edge of $uv$ labeled as $v_1$ and
$v_2$. Then $G$ has a 2-factor in which some cycle has $uv$ as a chord iff the dual $G'$ has a special
3-coloring.
\end{lem}

\noindent {\bf Proof:}
Suppose $G$ has such a 2-factor. We first assign colors $a,c$ as in the proof of Lemma~\ref{Pet} so that
no face in $G'$ is monochromatic. Let $C$ be the cycle in the 2-factor such that $uv$ is a chord of $C$. 
Without loss of generality, we can assume that the chord $uv$ is in the interior of $C$ and the region $R$ in 
the interior of $C$ that has $C$ as part of its boundary is colored $a$. In particular, this implies $v_1,v_2$ 
are colored $a$. The chord $uv$ splits $R$ into two regions, $R_1$ which contains the face corresponding to 
$v_1$ and $R_2$ which contains the face corresponding to $v_2$. We now change the color of all vertices in $G'$ 
corresponding to faces contained in the region $R_2$ to $b$. This gives a special 3-coloring of the dual 
triangulation $G'$. Conversely, if there exists  a special 3-coloring of $G'$, the edges in $G$ corresponding 
to the edges in $G'$ whose endpoints have the same color, and the edge $uv$, form a matching in $G$. Since 
there is no 2-edge-cut in $G$ containing the edge $uv$, there is a single edge between $v_1$ and $v_2$ in $G'$. 
This edge must be a bridge in the subgraph of $G'$ containing the edges having endpoints of the same color and 
the edge $v_1v_2$, since deleting it separates $v_1$ from $v_2$. Therefore the dual edge $uv$ must be a chord 
of some cycle, since it must become a self-loop after contracting all the edges in the 2-factor.
\hfill $\Box$

We will be working with plane triangulations and we define some simple operations on them. Let $uv$ be
an edge in a plane triangulation $G$. Let $w_1,w_2$ be the vertices such that $u,v,w_1$ and $u,v,w_2$
are the two faces in $G$ whose boundary contains the edge $uv$. Suppose $w_1 \neq w_2$. The operation of
\emph{flipping} the edge $uv$ replaces the edge $uv$ by the edge $w_1w_2$ on the same side of the
cycle $u,w_1,v,w_2$ that contains the edge $uv$. Note that the result of this operation is again
a triangulation. The operation of \emph{contracting} the edge $uv$ to the vertex $u$ deletes the edges
$vw_1$, $vw_2$ and $uv$ and merges the vertex $v$ with the vertex $u$. All the other edges incident with 
$v$ are now moved to $u$, preserving the circular order of the edges at all vertices and thus the embedding. 
If there was another edge between $u$ and $v$, it will now become a self-loop at $u$. The only difference
between ordinary contraction of an edge and this is that we delete the edges $vw_1$ and $vw_2$, so that
all faces in the resulting graph are still bounded by 3 edges. In the dual, this is equivalent to
deleting an edge and suppressing the resulting degree 2 vertices by contracting an edge incident with
them. Again, with this definition of contraction, the result of applying this operation is
a triangulation. The operations are also defined in the same way for internal edges of near-triangulations.

\section{Main Result}

In this section we prove the main theorem.

\begin{thm}
\label{main}
Every edge in a 2-edge-connected planar cubic graph is either contained in a 2-edge-cut or is a chord
of some cycle that is contained in a 2-factor of the graph.
\end{thm}

We will actually prove the following reformulation, whose equivalence follows from Lemma~\ref{3col}.

\begin{thm}
\label{Equiv}
Let $G$ be a plane triangulation with the boundary $v_1,v_2,v_3$. Then there exists a special
3-coloring of $G$ with the boundary assigned the colors $abc$.
\end{thm}

In order to prove  Theorem~\ref{Equiv}, we need to prove another statement, involving near-triangulations
with a 4-sided boundary, simultaneously by induction.

\begin{lem}
\label{4sides}
Let $G$ be a near-triangulation with boundary $v_1,v_2,v_3,v_4$. Then at least two of the assignments of
colors in $S = \{abbc, abca, abcc \}$ to the boundary are feasible for $G$. 
\end{lem}

We can classify 4-sided near-triangulations into 3 types depending on which two of the three
assignments in $S$ are feasible for them. Note that for some near-triangulations all 3 assignments may be
feasible, in which case it can be assumed to be any type. The three types correspond to the simplest 
4-sided near-triangulations, two with no internal vertex and a chord $v_1v_3$ or $v_2v_4$, and the other 
with exactly one internal vertex adjacent to all four vertices on the boundary. We show that every 4-sided 
near-triangulation is of one of these types. We will say a near-triangulation is of type $\mathcal{T}_1$ if 
the assignments $abbc, abca$ are feasible for it. It is of type $\mathcal{T}_2$ if 
the assignments $abca, abcc$ are feasible for it and of type $\mathcal{T}_3$ if $abbc, abcc$ are feasible
assignments.

\noindent {\bf Proof:} 
Suppose the theorem and/or the lemma is false and consider a counterexample $G$ to either with the minimum
number of vertices.

We first show that $G$ must be simple and cannot contain a separating triangle, that is, a triangle
that has vertices of $G$ in its interior as well as exterior.

Suppose $G$ has multiple edges, and let $e_1,e_2$ be two edges in $G$ with the same endpoints
$p,q$. Let $C$ be the 2-cycle containing the two edges. Since every face is 3 sided,
there must be vertices in the interior as well as exterior of $C$ and at least one of the edges, say $e_1$,
is an internal edge. Let $G'$ be the graph obtained from $G$ by deleting all the vertices and edges in the 
interior of $C$ and the edge $e_1$. Then $G'$ is also a triangulation or a near-triangulation with the same 
boundary but fewer vertices than $G$. We show that any assignment of colors to the boundary that is 
feasible for $G'$ is also feasible for $G$. Suppose $f'$ is a special 3-coloring of $G'$.
Let $G''$ be the triangulation obtained from $G$ by deleting the vertices and edges in the exterior of the 
cycle $C$ and the edge $e_1$, and let the boundary of $G''$ be the boundary of a face that contains $e_2$.
Relabel the vertex $p$ as $v_1$ and $q$ as $v_2$ in $G''$. Suppose $\{f'(p)\} \cup \{f'(q)\} \neq \{a,b\}$. 
Without loss of generality, assume $b \not\in \{f'(p)\} 
\cup \{f'(q)\}$. Lemma~\ref{Pet} implies, after swapping colors if necessary, that 
there exists a non-monochromatic 2-coloring $f''$ of $G''$ with colors $\{a,c\}$ such that $f''(p) = f'(p)$ and 
$f''(q)=f'(q)$. If $\{f'(p)\} \cup \{f'(q)\} = \{a,b\}$, then we may assume $p=v_1$ and $q=v_2$. Again by the 
minimality of $G$, $G''$ has a special 3-coloring $f''$ such that $f''(p) = f'(p)$ and $f''(q) = f'(q)$. In 
both cases, setting $f(v) = f'(v)$ for all $v \in V(G')$ and $f(v) = f''(v)$ for all $v \in V(G'')$ gives a 
special 3-coloring of $G$ with the same assignment of colors to the boundary as in $G'$. This contradicts the 
assumption that $G$ is a counterexample.

Suppose $G$ is simple but contains a separating triangle $C$ with vertices $p, q, r$.
Again, let $G'$ be the triangulation obtained from $G$ by deleting the vertices and edges in the interior of
$C$ and $G''$ the triangulation obtained from $G$ by deleting the vertices and edges in the exterior of $C$. 
Suppose $f'$ is any special 3-coloring of $G'$ with some assignment of colors to the boundary. Since $C$ is 
an internal face in $G'$, $|\{f'(p)\} \cup \{f'(q)\} \cup \{f'(r)\}| \ge 2$. If $|\{f'(p)\} \cup \{f'(q)\} \cup 
\{f'(r)\}| = 2$, then $\{f'(p)\} \cup \{f'(q)\} \cup \{f'(r)\} \neq \{a,b\}$, since we may assume $r \not\in 
\{v_1,v_2\}$ and either $f'(r) \neq f'(p)$ or $f'(r) \neq f'(q)$. Again, we may assume without loss of 
generality that the color $b$ is not assigned to any of the vertices $p,q,r$ in $f'$. Lemma~\ref{Pet} then 
implies, again after swapping colors and relabeling vertices if necessary, that there exists a 
non-monochromatic 2-coloring $f''$ of $G''$ with colors $\{a,c\}$ such that $f''(v) = 
f'(v)$ for all $v \in \{p,q,r\}$. If $|\{f'(p)\} \cup \{f'(q)\} \cup \{f'(r)\}| = 3$, we may assume without 
loss of generality that $p = v_1$, $q=v_2$, $f'(p) = a$, $f'(q) = b$ and $f'(r) = c$. Considering $v_1,v_2,r$ 
to be the boundary of $G''$, by the minimality of $G$, there is a special 3-coloring $f''$ of $G''$ such that 
$f''(v) = f'(v)$ for all $v \in \{p,q,r\}$. In either case, setting $f(v) = f'(v)$ for all $v \in V(G')$ and 
$f(v) = f''(v)$ for all $v \in V(G'')$ gives a special 3-coloring of $G$ with the same assignment of colors to 
the boundary as in $G'$. Again, this contradicts the assumption that $G$ is a counterexample. 

We now assume $G$ is simple and has no separating triangles and consider cases depending on whether
$G$ is a triangulation or a near-triangulation with a 4-sided boundary.

\noindent{\bf Case 1:} Suppose $G$ is a triangulation. If $G$ contains a vertex of degree 2, since 
there are no multiple edges, $G$ must be $K_3$, in which case, the assignment $abc$ of colors to the boundary 
is a special 3-coloring of $G$. If $G$ contains a vertex of degree 3, since there are no separating triangles, 
$G$ must be $K_4$. However, in this case, assigning color $c$ to the internal vertex and $abc$ to the boundary 
gives a special 3-coloring of $G$. So we may assume every vertex in $G$ has degree at least 4. 

We next show that $v_3$ must have degree exactly 4. Suppose not and let $v_1, v_4, v_5, v_6$ be the
neighbors of $v_3$ in circular order, such that $v_1,v_3,v_4$, and $v_3,v_4,v_5$ are 
internal faces in $G$. In any special 3-coloring of $G$ with $abc$ assigned to the boundary, $v_3$ must 
be assigned color $c$. The vertex $v_5$ is not adjacent to either $v_1$ or $v_2$, otherwise there
is a separating triangle in $G$, separating $v_4$ and $v_6$.  Let $G'$ be the triangulation obtained from 
$G$ by flipping the edge $v_3v_4$, replacing it by $v_1v_5$, and then contracting $v_1v_5$ to the 
vertex $v_1$. Then $G'$ has the same boundary as $G$ and no self-loops, hence by the minimality of $G$, it 
has a special 3-coloring $f'$ with the colors $abc$ assigned to the boundary. Define $f(v_5) = f'(v_1) = a$ 
and $f(v) = f'(v)$ for all $v \in V(G')$.  We show that $f$ is a special 3-coloring of $G$. All faces in $G$ 
except $v_1,v_3,v_4$ and $v_3,v_4,v_5$ correspond to faces in $G'$, with vertex $v_5$ replaced by $v_1$ if it
appears on the boundary of the face. Since $f(v_5) = f'(v_1)$ and $f(v) = f'(v)$ for all other vertices,
these faces are non-monochromatic in $f$. Since $f(v_3) = f'(v_3) = c$ and $f(v_5) = f(v_1) = f'(v_1) = a$, the 
other two faces are also non-monochromatic. The only pairs of adjacent vertices in $G$ that are not adjacent 
in $G'$ are $v_3,v_4$ and those that include $v_5$ and one of its neighbors. Since $v_3$ is colored $c$, it 
cannot form an $ab$-pair with any other vertex. Any neighbor of $v_5$ in $G$ is a neighbor of $v_1$ in $G'$. 
The only neighbor of $v_1$ that is colored $b$ is $v_2$, which is not a neighbor of $v_5$. Therefore there are 
no $ab$-pairs in $f$ apart from $v_1,v_2$ and $f$ is a special 3-coloring of $G$, a contradiction. Note that 
this argument fails if the degree of $v_3$ is exactly 4, since $v_5$ is now adjacent to $v_2$. 

We now assume the degree of $v_3$ is exactly 4 and as before, let $v_1,v_4,v_5,v_2$ be its
four neighbors in circular order. Let $G'$ be the near-triangulation obtained from $G$ by deleting the vertex 
$v_3$, having the boundary $v_1, v_2, v_5, v_4$ in circular order. The minimality of $G$ implies $G'$ is
of one of types $\mathcal{T}_1$, $\mathcal{T}_2$, or $\mathcal{T}_3$. In all 3 cases, there exists
a special 3-coloring $f'$ of $G'$ with the boundary assigned colors $abbc$ or $abca$. In either case,
setting $f(v_3) = c$ and $f(v) = f'(v)$ for all $v \in V(G')$ gives a special 3-coloring of $G$, with
colors $abc$ assigned to the boundary. This contradicts the assumption that $G$ is a counterexample.

\noindent {\bf Case 2:} Suppose $G$ is a 4-sided near-triangulation.
Suppose $G$ has a chord, in which case there are no internal vertices in $G$. If the chord is $v_1v_3$,
the assignments $abca$ and $abcc$ are special 3-colorings of $G$, but $abbc$ is not since the vertex
$v_1$ with color $a$ is adjacent to the vertex $v_3$ of color $b$. In this case, $G$ is of type
$\mathcal{T}_2$. If the chord is $v_2v_4$ then the assignments $abbc$ and $abcc$ are special 
3-colorings and $G$ is of type $\mathcal{T}_3$.

We may assume $G$ has no chords. Let $v_5$ be the internal vertex in $G$ such that $v_1,v_2,v_5$ is the
internal face in $G$ containing the edge $v_1v_2$. If $v_5$ is adjacent to both $v_3$ and $v_4$, then
there are no other vertices in $G$. Assigning color $c$ to $v_5$ and either colors $abbc$ or $abca$
to the boundary, gives a special 3-coloring of $G$. Thus $G$ is of type $\mathcal{T}_1$.

Suppose $v_5$ is adjacent to $v_3$ but not to $v_4$. Let $G'$ be the triangulation obtained from $G$
by deleting the vertex $v_2$ and adding the edge $v_1v_3$ in the exterior of the cycle $v_1,v_5,v_3,v_4$.
Then $G'$ is a triangulation with boundary $v_1,v_3,v_5$ and by the minimality of $G$, $G'$ has a
special 3-coloring $f'$ with colors $abc$ assigned to the boundary. Since $v_4$ is adjacent to both
$v_1$ and $v_3$, we must have $f'(v_4) = c$. Setting $f(v_2) = b$ and $f(v) = f'(v)$ for all $v \in V(G')$
gives a special 3-coloring of $G$ with colors $abbc$ assigned to the boundary. Similarly, Lemma~\ref{Pet}
implies that $G'$ has a non-monochromatic 2-coloring $f''$ with colors $acc$ assigned to the boundary.
In this case, $v_4$ may be colored either $a$ or $c$. Setting $f(v_2) = b$ and $f(v) = f''(v)$ for all
$v \in V(G')$, gives a special 3-coloring of $G$ with either $abca$ or $abcc$ assigned to the boundary.
Then $G$ is of type $\mathcal{T}_1$ or $\mathcal{T}_3$, depending on the color of $v_4$ in $f''$.
A symmetrical argument holds if $v_5$ is adjacent to $v_4$ but not to $v_3$.

We may now assume $v_5$ is not adjacent to either $v_3$ or $v_4$. The degree of $v_5$ is at least four,
since $G$ has no separating triangles. Again, we show that the degree must be exactly four, using
essentially the same argument as in the case of triangulations. Note that in any special 3-coloring
of $G$, $v_5$ must  be assigned color $c$. Suppose the degree of $v_5$ is at least five and let
$v_2,v_6,v_7,v_8$ be vertices adjacent to $v_5$ in circular order, such that $v_2,v_5,v_6$ and $v_5,v_6,v_7$ 
are internal faces in $G$. The vertex $v_7$ cannot be adjacent to either $v_1$ or $v_2$, otherwise $G$ has a 
separating triangle. Let $G'$ be obtained from $G$ by flipping the edge $v_5v_6$, replacing it by $v_2v_7$ and 
then contracting the edge $v_2v_7$ to the vertex $v_2$. Then $G'$ has the same boundary as $G$ and fewer
vertices, and by the minimality of $G$, is of one of the three types. As argued in the case of
triangulations, any special 3-coloring of $G'$ with an assignment of colors from $S$ to the boundary,
can be extended to a special 3-coloring of $G$ with the same assignment to the boundary,
by assigning the color of $v_2$ ($b$) to $v_7$. Therefore $G$ is of the same type as $G'$, a contradiction.

Suppose the degree of $v_5$ is exactly 4 and let $v_1,v_2,v_6,v_7$ be its neighbors in circular order.
The previous argument cannot be used now as $v_7$ is adjacent to $v_1$ and $v_1,v_7$ would become
an $ab$-pair in the coloring $f$ of $G$.

Suppose there exists a cycle $C$ of length four $v_1,v_2,p,q$ in $G$ such that $\{p,q\} \not\subseteq 
\{v_5,v_6,v_7\}$ and $\{p,q\} \neq \{v_3,v_4\}$. Let $G'$ be the near-triangulation obtained from
$G$ by deleting the vertices in the exterior of the cycle $C$, with the boundary $v_1,v_2,p,q$. Since at 
least one of $v_3,v_4$ must be in the exterior of $C$, $G'$ has fewer vertices than $G$. The minimality of 
$G$ implies $G'$ is of one of the three types. We replace the vertices in the interior of $C$ by an equivalent 
smaller subgraph of the same type. Note that at least two of the vertices $v_5,v_6,v_7$ must be contained
in the interior of $C$. Let $H$ be the graph obtained from $G$ by deleting the vertices in the
interior of $C$. If $G'$ is of type $\mathcal{T}_1$, add a new vertex $r$ to $H$ in the interior of
the cycle $C$ and join it to all four vertices $v_1,v_2,p,q$ and call the resulting graph $G''$.
If $G'$ is of type $\mathcal{T}_2$, add the edge $v_1p$ in the interior of $C$  to construct $G''$ from
$H$. If $G'$ is of type $\mathcal{T}_3$, add the edge $v_2q$ in the interior of $C$ to construct $G''$
from $H$. In all cases, $G''$ is a near-triangulation with the same boundary but fewer vertices than
$G$, and by the minimality of $G$, is of one of the three types. We claim that $G$ is of the same type
as $G''$. Let $f''$ be any special 3-coloring of $G''$. If $G''$ was constructed by adding the vertex
$r$ in the interior of $C$,  since $r$ is adjacent to $v_1$ and $v_2$, we must have $f''(r) = c$.
This implies that either $f''(p) = b$ and $f''(q) = c$ or $f''(p) = c$ and $f''(q) = a$. Since in
this case $G'$ was of type $\mathcal{T}_1$, there exists a special 3-coloring $f'$ of $G'$ such that
$f'(p) = f''(p)$ and $f'(q) = f''(q)$. Similarly, in the other cases, it can be argued that the
possible values of $f''(p)$ and $f''(q)$ are such that there exists a special 3-coloring $f'$ of
$G'$ such that $f'(p) = f''(p)$ and $f'(q) = f''(q)$. In all cases, we have $f'(v_1) = f''(v_1) = a$
and $f'(v_2) = f''(v_2) = b$. Defining $f(v) = f'(v)$ for all $v \in V(G')$ and $f(v) = f''(v)$ for
all $v \in V(G'')$ gives a special 3-coloring of $G$ with the same assignment of colors to the
boundary as in the coloring of $G''$. This implies $G$ is of the same type as $G''$, a contradiction.

We may now suppose there is no such cycle $C$ in $G$. This implies that $v_7$ is not adjacent to
$v_3$ and $v_6$ is not adjacent to $v_4$. Suppose $v_7$ is adjacent to $v_4$ and $v_6$ is adjacent to
$v_3$. We claim that $G$ is of type $\mathcal{T}_1$ in this case. Let $G'$ be the triangulation
obtained from $G$ by deleting the vertices $v_1, v_2, v_5$ and adding the edge $v_3v_7$ in the
exterior of the cycle $v_3,v_4,v_7,v_6$. Let the boundary of $G'$ be $v_7,v_3,v_4$. The minimality
of $G$ implies that $G'$ has a special 3-coloring $f'$ with colors $abc$ assigned to the boundary.
Since $v_6$ is adjacent to both $v_3$ and $v_7$, we have $f'(v_6) = c$. Setting $f(v_1) = a$,
$f(v_2) = b$, $f(v_5) = c$ and $f(v) = f'(v)$ for all $v \in V(G')$ gives a special 3-coloring of
$G$ with colors $abbc$ assigned to the boundary. Similarly, by deleting the vertices $v_1,v_2,v_5$ and adding 
the edge $v_4v_6$ to construct a triangulation $G''$ with boundary $v_4,v_6,v_3$, we can find a special 
3-coloring of $G$ with colors $abca$ assigned to the boundary. Thus $G$ is of type $\mathcal{T}_1$,
a contradiction.

Suppose, without loss of generality, that $v_6$ is not adjacent to $v_3$. The argument in the other
case is symmetrical, after relabeling the vertices and swapping the colors $a$ and $b$. Let $w \neq v_5$
be the internal vertex in $G$ such that $v_2,v_6,w$ is the other face in $G$ containing the edge 
$v_2v_6$. We show that we can find special 3-colorings of $G$ with the additional condition that
$w$ is colored $c$. Let $v_2 = w_1$, $w_2, \ldots, w_d = v_6$ be the vertices adjacent to $w$ in circular
order, starting with $v_2$ and ending with $v_6$, where $d \ge 4$ is the degree of $w$. Note that
if $w$ is adjacent to $v_3$, we must have $w_2 = v_3$, otherwise $G$ has a separating triangle.
The vertex $v_1$ cannot be adjacent to any vertex $w_i$ for $i > 1$, otherwise $v_1,v_2,w,w_i$ is a 4-cycle
in $G$. This implies $v_4, v_7$ are not adjacent to $w$. Also, no $w_i$ can be adjacent to a $w_j$
for $1 < |i-j| < d-1$, otherwise $G$ has a separating triangle. 

Construct a near-triangulation $G'$ from $G$ as follows. First flip all the edges $ww_{2i}$ and replace
them by the edges $w_{2i-1}w_{2i+1}$ for $1 \le i < \lceil d/2 \rceil$. Then contract all the
edges $w_{2i-1}w_{2i+1}$ for $1 \le i < \lceil d/2 \rceil$ to the vertex $w_1 = v_2$. 
If $d$ is even, the vertex $w$ will have degree 2 in the resulting graph with $v_2$ and $v_6$ as its two 
neighbors. There will be two edges between $v_6$ and $v_2$, corresponding to the edges $v_2v_6 = w_1w_d$ and 
$w_{d-1}w_d$ in $G$, and $w$ is the only vertex contained in the interior of the 2-cycle formed by the two 
edges. There will be no self-loops in the resulting graph. If $d$ is odd, $w$ will have degree 1 with $v_2$ 
as its only neighbor. There will be a self-loop at $v_2$ corresponding to the edge $w_1w_d$ in $G$ whose
interior will contain the vertex $w$ and the edge $v_2w$. The other face containing this loop will contain two 
edges between $v_2$ and $v_5$, corresponding to the original edges $v_2v_5$ and $v_5v_6$ in $G$. In this case, 
we delete the vertex $w$, the self loop at $v_2$ and one of the two edges between $v_2$ and $v_5$ that are 
in the same face as the self-loop at $v_2$. Let $G'$ be the resulting graph. 

Note that $G'$ has the same boundary as $G$ as neither $v_3$ nor $v_4$ will be merged with $v_2$. There
are no self-loops in $G'$ and since at least one edge will be contracted, $G'$ has fewer vertices than $G$.  
The minimality of $G$, implies that $G'$ is of one of the three types. We show that $G$ is of the same type as 
$G'$. Consider any special 3-coloring $f'$ of $G'$ with some assignment of colors in $\{abbc, abca, abcc \}$ to 
the boundary. In the case $d$ is even, since the two faces whose boundary contains $w$ both have $v_2$ and 
$v_6$ as the other two vertices, and $f'(v_2) = b$, we can set $f'(w) = c$, irrespective of the color of $v_6$. 
Now we claim that setting $f(w_{2i+1}) = b$ for $1 \le i < \lceil d/2 \rceil$, $f(w) = c$ and $f(v) = f'(v)$ 
for all $v \in V(G')$ gives
a special 3-coloring of $G$ with the same assignment of colors to the boundary as in $f'$. Any face in $G$
that does not include $w$ on its boundary, except possibly the face $v_2,v_5,v_6$ when $d$ is odd, 
corresponds to a face in $G'$ with the vertex $w_{2i+1}$ replaced by $w_1 = v_2$ if $w_{2i+1}$ is on its
boundary for some $1 \le i < \lceil d/2 \rceil$. Since all vertices $w_{2i+1}$ have the same color $b$ as 
$w_1$, these are non-monochromatic in $f$ as well. The faces in $G$ that include $w$ are non-monochromatic 
since $f(w) = c$ and $f(w_{2i-1}) = b$ for $1 \le i \le \lceil d/2 \rceil$. The face $v_2,v_5,v_6$ is 
non-monochromatic since $f(v_5)$ must be $c$. There cannot be any $ab$-pairs in $f$ apart from $v_1,v_2$, since 
the only vertex of color $a$ adjacent to $v_2$ in $G'$ is $v_1$, but none of the vertices $w_{2i+1}$ for 
$1 \le i < \lceil d/2 \rceil$ is adjacent to $v_1$. Thus $G$ has a special 3-coloring with the same assignment
of colors to the boundary as in $G'$ and hence is of the same type as $G'$.

This completes the proof.
\hfill $\Box$

\section{Disconnected 2-factors}

In this section, we consider disconnected 2-factors in 2-edge-connected planar cubic graphs. This is
equivalent to finding separating perfect matchings in the graph. While the existence of these in all 
2-edge-connected planar cubic graphs with at least six vertices was shown in~\cite{D}, here we consider their 
existence with the additional restriction that they should include a specified edge. We give a constructive 
characterization of planar cubic 2-edge-connected graphs that contain an edge such that every 2-factor 
containing the edge is a Hamiltonian cycle. We also characterize 2-edge-connected planar cubic graphs that
have an edge that is not contained in a separating perfect matching.

We first show using Theorem~\ref{main} that cyclically 4-edge-connected planar cubic graphs except $K_2^3$ and 
$K_4$ do not contain any such edges.

\begin{thm}
\label{disconn}
Let $G$ be a 2-edge-connected planar cubic graph and let $uv$ be any edge in $G$ such that there is no edge
parallel to it and $G-\{u,v\}$ is 2-edge-connected. Then there exists a separating perfect matching in $G$ 
including the edge $uv$. 
\end{thm}

\noindent {\bf Proof:}
Consider any plane embedding of $G$. Since there is no edge parallel to $uv$, both $u,v$ must have 2 neighbors
each not in $\{u,v\}$. Let $u_1,u_2$ be the other neighbors of $u$ and $v_1,v_2$ the other neighbors of $v$ 
such that $u_1,u,v,v_1$ are consecutive vertices on the boundary of a face of $G$, as are $u_2,u,v,v_2$. The 
vertices $u_1,u_2,v_1,v_2$ must all be distinct otherwise there is a bridge in $G-\{u,v\}$. Let $G'$ be the 
planar cubic graph obtained from $G-\{u,v\}$ by adding two new vertices $u',v'$ adjacent to each other, joining 
$u'$ to $u_1,v_1$ and $v'$ to $u_2,v_2$. Since $G-\{u,v\}$ is 2-edge-connected, $u'v'$ is not contained in a 
2-edge-cut in $G'$. Theorem~\ref{main} implies there is a 2-factor in $G'$ such that $u'v'$ is a chord of some 
cycle $C$ in the 2-factor. The cycle $C$ must contain the edges $u_1u'$, $u'v_1$, $u_2v'$ and $v'v_2$. Then 
planarity implies $C -\{u',v'\}$ is the disjoint union of a path from $u_1$ to $u_2$ and a path from $v_1$ to 
$v_2$. Adding the edges $u_1u$, $u_2u$, $v_1v$ and $v_2v$ to these two paths gives two disjoint cycles in $G$ 
which include all vertices in $C$ except $u',v'$. These two cycles, along with the cycles other than $C$ in the 
2-factor of $G'$ (if any), give a disconnected 2-factor in $G$ not containing the edge $uv$. The complement of 
this is a separating perfect matching in $G$ containing $uv$.
\hfill $\Box$

An immediate consequence of Theorem~\ref{disconn} is that every edge in a cyclically 4-edge-connected
planar cubic graph, other than $K_2^3$ and $K_4$, is contained in a separating perfect matching. Any edge 
$uv$ in such a graph satisfies the condition $G-\{u,v\}$ is 2-edge-connected. This also implies that every
edge in such a graph is contained in a disconnected 2-factor. The result in~\cite{D} that every 
2-edge-connected planar cubic graph except $K_2^3$ and $K_4$ has a disconnected 2-factor also follows
easily from this. If the graph is not cyclically 4-edge-connected, it must have a cyclic cut of size at most 3.
If it has a 2-edge-cut then any perfect matching containing an edge in the cut must also contain the other
edge and is a separating perfect matching. Petersen's theorem implies there exists such a matching. On
the other hand, if the graph $G$ is 3-edge-connected, any edge $uv$ contained in a cyclic cut of size 3
must satisfy $G-\{u,v\}$ is 2-edge-connected, otherwise there is a 2-edge-cut in $G$.

We next consider 3-edge-connected planar cubic graphs and characterize those that have an edge
that is not contained in a disconnected 2-factor and also those that have an edge not contained in a
separating perfect matching. Let $H_n$ for $n \ge 1$ denote the graph obtained from a cycle
$v_0,v_1,\ldots,v_{2n-1}$ of length $2n$ by adding the edges $v_0v_n$ and $v_iv_{2n-i}$ for $1 \le i < n$.
Note that $H_2$ is $K_4$ and $H_3$ is the prism $K_3 \times K_2$.

\begin{thm}
\label{3conn}
Let $G$ be a 3-edge-connected planar cubic graph with an edge $uv$ such that every 2-factor of $G$ that
contains $uv$ is a Hamiltonian cycle. Then $G$ is $H_n$ for some $n \ge 1$ and $uv$ is the edge $v_0v_n$.
\end{thm}

\noindent {\bf Proof:}
The proof is by induction on the number of vertices. If $G$ has only two vertices, it must be $K_2^3$ and the 
result holds. Suppose $G$ has $2n$ vertices for some $n > 1$. Let $u_1,u_2$ be the neighbors of $u$ other than 
$v$. Theorem~\ref{disconn} implies $G-\{u,u_1\}$ as well as $G-\{u,u_2\}$ 
contains a bridge, otherwise there is a separating perfect matching containing $uu_1$ or $uu_2$, whose 
complement is a disconnected 2-factor in $G$ containing the edge $uv$, a contradiction. Let $pq$ be a bridge in 
$G-\{u,u_1\}$, let $C_1$, $C_2$ be the components of $(G-\{u,u_1\})-pq$ and assume without loss of 
generality that $p$ is in $C_1$ and $q$ in $C_2$. Since $G$ is 3-edge-connected, $u_2$ and $v$
belong to different components, and without loss of generality, $u_2$ is in $C_1$. Note that $u_1$ also has
one neighbor in $C_1$ and one in $C_2$. If the component $C_1$ is non-trivial, it must contain at least
3 vertices and must be 2-connected, otherwise $G$ has a 2-edge-cut. The same holds for $C_2$. This implies
that if $C_1$ is not trivial, we can find a 2-edge-connected subgraph in $G-\{u,u_2\}$ containing both
$v$ and $u_1$. This is obtained by adding a path to $C_2$ that goes through $u_1$ and $pq$, and uses a
path in $C_1-\{u_2\}$ between $p$ and the vertex in $C_1$ that is adjacent to
$u_1$. This contradicts the fact that there must be a bridge that separates $v$ and $u_1$ in $G-\{u,u_2\}$.
Therefore $C_1$ must be trivial and contain only the vertex $p = u_2$. This implies $u,u_1,u_2$ form a
triangle in $G$. Let $G'$ be obtained from $G$ by contracting the triangle $u,u_1,u_2$ to the vertex $u$.
If $uv$ is a contained in a disconnected 2-factor in $G'$, then it is also contained in a disconnected
2-factor in $G$, obtained by adding appropriate edges from the contracted triangle. Otherwise $G'$ must
be $H_{n-1}$ and $uv$ must be the edge $v_0v_{n-1}$. Then the planarity of $G$ implies that 
$G$ must be $H_n$ and $uv$ the edge $v_0v_n$.
\hfill $\Box$

Theorem~\ref{disconn} implies that if $uv$ is an edge in a 3-edge-connected planar cubic graph other than
$K_2^3$ that is not contained in a separating perfect matching, then $G-\{u,v\}$ has a bridge say $pq$. 
Let $C_1,C_2$ be the components of $(G-\{u,v\})-pq$ and let $G_1$ and $G_2$ be the graphs obtained from $G$ 
by contracting the vertices $V(C_2) \cup \{u,v\}$ and $V(C_1) \cup \{u,v\}$ to the vertices $x$ and $y$, 
respectively. If there is a separating perfect matching in either $G_1$ containing the edge $px$ or one in 
$G_2$ containing the edge $qy$, then we get a separating perfect matching in $G$ including $uv$, a 
contradiction. Therefore both $G_1$ and $G_2$ are 3-edge-connected planar cubic graphs and the edges $px$ and 
$qy$ are not contained in a separating perfect matching in $G_1$ and $G_2$, respectively. The converse of this 
is also true and any graph $G$ constructed from $G_1$ and $G_2$ in this way has the edge $uv$ not contained 
in a separating perfect matching. The base case is if $G$ is $K_2^3$ and any other 3-edge-connected planar 
cubic graph with an edge not contained in a separating perfect matching can be obtained from $K_2^3$ using 
this operation.

Finally, we consider the 2-edge-connected case. Suppose $G$ is a 2-edge-connected planar cubic graph
with an edge $uv$ that is not contained in a disconnected 2-factor. If $uv$ is contained in a 2-edge-cut,
let $pq$ be the other edge in the cut. Let $C_1,C_2$ be the components of $G-\{uv, pq\}$ and without loss of
generality, assume $u,p \in V(C_1)$ and $v,q \in V(C_2)$. Let $G_1,G_2$ be the graphs obtained from $C_1, C_2$ 
by adding the edges $up, vq$, respectively. If there is a disconnected 2-factor in either $G_1$ that contains
the edge $up$, or a disconnected 2-factor in $G_2$ that contains $vq$, then we get a disconnected
2-factor in $G$ containing $uv$. Therefore both $G_1, G_2$ are 2-edge-connected planar cubic graphs
with edges $up, vq$ that are not contained in a disconnected 2-factor, respectively. The converse
of this also follows easily, and any graph $G$ constructed in this way, has the edge $uv$ not
contained in a disconnected 2-factor. 

Suppose the edge $uv$ is not contained in a 2-edge-cut. If $G$ is 3-edge-connected, it follows from 
Theorem~\ref{3conn} that $G$ is $H_n$ for some $n \ge 1$ and $uv$ is the edge $v_0v_n$.  Suppose there exists 
a 2-edge-cut $\{p_1q_1, p_2q_2\}$ in $G$ such that $p_1,p_2$ are vertices in the component $C_1$ of 
$G -\{p_1q_1,p_2q_2\}$ that contains the edge $uv$. Let $C_2$ be the other component. Choose such a cut
such that the size of $C_2$ is as large as possible. Let $G_1,G_2$ be the graphs obtained from $C_1,C_2$ 
by adding the edges $p_1p_2, q_1q_2$, respectively. Then $uv$ is not contained in a 2-edge-cut in $G_1$, 
otherwise we get a 2-edge-cut in $G$ containing $uv$. Also, if there is a 2-edge-cut in $G_1$, the vertices
$u,v,p_1,p_2$ must be contained in the same component of the cut, otherwise it contradicts the assumption
that the size of $C_2$ was as large as possible. Repeat this process of 2-edge-cut reduction on $G_1$ until we 
are left with a 3-edge-connected graph $G_1$ that contains the edge $uv$. Thus the graph $G$ may be viewed as
being obtained from a 3-edge-connected planar cubic graph $G_1$ by replacing some of the edges by the graphs
$G_2$. 

If the final graph $G_1$ has a disconnected 2-factor containing the edge $uv$, we can find one in $G$ as well.
If $p_1p_2$ is any added edge in the graph $G_1$, if the 2-factor does not contain the edge $p_1p_2$, adding
a 2-factor in $G_2$ not containing the edge $q_1q_2$ gives a disconnected 2-factor in $G$. If $p_1p_2$ is
included in the 2-factor, take a 2-factor in $G_2$ containing the edge $q_1q_2$. Replacing the edges
$p_1p_2, q_1q_2$ by the edges $p_1q_1$, $p_2q_2$ gives a disconnected 2-factor in $G$. Theorem~\ref{3conn}
then implies the final graph $G_1$ must be $H_n$ for some $n \ge 1$. The same argument holds if any of the 
graphs $G_2$ that replace an edge in $G_1$ has a disconnected 2-factor containing the edge $q_1q_2$. Take any 
2-factor in $G_1$ containing the edge $uv$. If the edge $p_1p_2$ is not in the 2-factor, we add any 2-factor 
in $G_2$ not containing $q_1q_2$, otherwise add a disconnected 2-factor containing $q_1q_2$ and replace the 
edges $p_1p_2, q_1q_2$ by $p_1q_1, p_2q_2$. Finally, if $G_1$ has a 2-factor containing $uv$ but not 
containing any one of the added edges $p_1p_2$, we get a disconnected 2-factor by adding a 2-factor in $G_2$ 
not containing $q_1q_2$. Therefore the only case when we cannot get a disconnected 2-factor is if $G_1$ is
$H_n$ for some $n \ge 1$, the edge $uv$ is $v_0v_n$ and any 2-factor in $H_n$ that contains $v_0v_n$
also includes the replaced edges $p_1p_2$. There are only two perfect matchings in $H_n$ that do
not contain the edge $v_0v_n$, depending on whether it includes $v_0v_1$ or $v_0v_{2n-1}$. In either
case, the matchings cannot include the edges $v_iv_{2n-i}$ for $1 \le i < n$ and any 2-factor that
includes $v_0v_n$ must include these edges. These are precisely the edges that can be replaced by
the graphs $G_2$, and $G$ must have been obtained from some $H_n$ for $n \ge 1$ by these operations. The
edge $uv$ is still the edge $v_0v_n$ of $H_n$. 

This describes the structure of all 2-edge-connected planar cubic graphs that have an edge $uv$ such that every 
2-factor containing the edge is a Hamiltonian cycle. We now consider the same for all such graphs that have 
an edge $uv$ that is not contained in any separating perfect matching.

Let $G$ be a 2-edge-connected planar cubic graph and $uv$ an edge in $G$ that is not contained
in any separating perfect matching. Suppose there exists another edge between $u$ and $v$. Then either
$G$ is $K_2^3$ or $u,v$ have distinct neighbors $u_1,v_1$, respectively, not in $\{u,v\}$.
Assume $G$ is not $K_2^3$ and let $G'$ be obtained from $G-\{u,v\}$ by adding the edge $u_1v_1$. Then
$G'$ is also a 2-edge-connected planar cubic graph and it is easy to see that $G$ has a separating perfect
matching including $uv$ iff $G'$ has a disconnected 2-factor containing the added edge $u_1v_1$. This gives
a way of constructing graphs having an edge not contained in a separating perfect matching from graphs
with an edge not contained in a disconnected 2-factor.

If $uv$ is contained in a 2-edge-cut then any perfect matching containing $uv$ must contain the
other edge in the 2-edge-cut, and any such perfect matching is a separating perfect matching. Petersen's
theorem implies there exists such a matching, a contradiction.

Suppose $G-\{u,v\}$ is disconnected and let $C_1,C_2$ be the components of $G-\{u,v\}$.
Let $u_1,u_2$ be the neighbors of $u$ and $v_1,v_2$ the neighbors of $v$ in $C_1,C_2$, respectively.
Let $G_1,G_2$ be the graphs obtained from $C_1,C_2$ by adding the edges $u_1v_1$, $u_2v_2$,
respectively. Now $G$ has a separating perfect matching containing the edge $uv$ iff either $G_1$ or
$G_2$ has a separating perfect matching not including $u_1v_1$ or $u_2v_2$, respectively. In other words,
there is no disconnected 2-factor in $G_1$ or $G_2$ containing the edges $u_1v_1$ or $u_2v_2$, 
respectively. The converse also follows similarly. This gives another way of constructing a graph having 
an edge not contained in a separating perfect matching from 2 smaller graphs having an edge not contained
in a disconnected 2-factor. 

If $G-\{u,v\}$ is connected but has a bridge, it follows from the same argument as in the
3-edge-connected case that $G$ can be constructed from two smaller graphs having an edge
not contained in a separating perfect matching. Finally, if $G-\{u,v\}$ is 2-edge-connected, since
there is no edge parallel to $uv$, it follows from Theorem~\ref{disconn} that there exists a
separating perfect matching including $uv$, a contradiction.

This covers all possibilities and gives a constructive characterization of 2-edge-connected planar
cubic graphs that have an edge not contained in a separating perfect matching.

\section{Extension}

It seems possible that the results in this paper can be extended further. In particular, it
appears to be true that for any two edges in a 2-edge-connected planar cubic graph, there
exists a 2-factor in the graph such that the two edges are contained in the same cycle of the 2-factor.
If the two edges are on the boundary of the same face, we can subdivide the two edges and add an
edge joining the two degree 2 vertices in the interior of the face. The added edge is not contained in
a 2-edge-cut in the resulting 2-edge-connected planar cubic graph, and Theorem~\ref{main} implies it is 
a chord of some cycle in a 2-factor of the graph. This gives a 2-factor in the original graph such
that both edges are contained in the same cycle. It is also possible that the same property holds for all 
connected cubic bipartite graphs. Here we prove it for connected planar cubic bipartite graphs.

\begin{thm}
Let $G$ be a connected planar cubic bipartite graph and let $e_1,e_2$ be any two edges in $G$. There
exists a 2-factor in $G$ such that both the edges are contained in the same cycle of the 2-factor.
\end{thm}

\noindent {\bf Proof:}
We prove by induction on the number of vertices. Note that since $G$ is cubic and bipartite, it cannot have
a bridge. Suppose $G$ has a 2-edge-cut $\{p_1q_1,p_2q_2\}$.
Let $C_1,C_2$ be the components of $G-\{p_1q_1,p_2q_2\}$ and assume $p_1,p_2 \in V(C_1)$ and $q_1,q_2
\in V(C_2)$. Let $G_1$, $G_2$ be the graphs obtained by adding the edges $p_1p_2$, $q_1q_2$ to $C_1,C_2$
respectively. We will call the edge $p_1p_2$ as $e_1$ if $e_1$ is in the 2-edge-cut and similarly
call $q_1q_2$ as $e_2$ if $e_2$ is in the 2-edge-cut. If $e_1,e_2$ are now contained in the same graph,
say $G_1$, by induction, there exists a 2-factor $F_1$ in $G_1$ such that $e_1,e_2$ are in the same cycle
of $F_1$. If $F_1$ does not contain the edge $p_1p_2$, let $F_2$ be any 2-factor in $G_2$ not containing
the edge $q_1q_2$. Then $F_1 \cup F_2$ is a 2-factor in $G$ such that $e_1$ and $e_2$ are contained in
the same cycle of the 2-factor. If $F_1$ contains $p_1p_2$ then let $F_2$ be a 2-factor in $G_2$
containing $q_1q_2$. Then $F = (F_1 \setminus \{p_1p_2\}) \cup (F_2 \setminus \{q_1q_2\}) \cup 
\{p_1q_1,p_2q_2\}$ is a 2-factor in $G$ such that $e_1,e_2$ are contained in the same cycle of $F$.
Finally, suppose $e_1,e_2$ are in different graphs, say $e_1$ is in $G_1$ and $e_2$ in $G_2$.
By induction, there is a 2-factor $F_1$ in $G_1$ such that $e_1$ and $p_1p_2$ are in the same cycle of
$F_1$. Similarly, there is a 2-factor $F_2$ in $G_2$ such that $e_2$ and $q_1q_2$ are in the same cycle
of $F_2$. Then $F = (F_1 \setminus \{p_1p_2\}) \cup (F_2 \setminus \{q_1q_2\}) \cup \{p_1q_1,p_2q_2\}$ is a 
2-factor in $G$ such that $e_1,e_2$ are contained in the same cycle of $F$.

We may now assume $G$ is 3-edge-connected. If $G$ is $K_2^3$ the result is obviously true. If $G$ has at
least 4 vertices then there are no multiple edges and since $G$ is bipartite, every face in a plane embedding 
of $G$ has even size at least four. Euler's formula
implies $G$ has at least 6 faces of size 4, and hence at least two such that they do not contain
any of the edges $e_1,e_2$ on their boundary. Let $v_1,v_2,v_3,v_4$ be the vertices on the boundary
of such a face and let $u_i$ be the neighbor of $v_i$ not on the boundary of the face, for $1 \le i \le 4$.
The $u_i$ must be distinct vertices since $G$ is planar and bipartite. If $e_1$ is one of the edges $u_iv_i$, we may
assume without loss of generality, it is $u_1v_1$. If $e_2$ is also incident with some vertex $v_i$,
we may assume it is one of $u_2v_2$ or $u_3v_3$. Let $G'$ be the graph obtained from $G$ by deleting
the vertices $v_1,v_2,v_3,v_4$ and adding the edges $u_1u_4$ and $u_2u_3$. Then $G'$ is a cubic bipartite
graph and must be connected, otherwise $G$ has a 2-edge-cut. Label the added edge $u_1u_4$ as $e_1$ if
$e_1$ was the edge $u_1v_1$. Similarly, label $u_2u_3$ as $e_2$ if either $u_2v_2$ or $u_3v_3$ was $e_2$.
By induction, $G'$ has a 2-factor $F'$ such that $e_1$ and $e_2$ are in the same cycle of $F'$. If $F'$
does not contain any of the added edges, then adding the edges in the 4-cycle $v_1,v_2,v_3,v_4$ to $F'$ gives
a 2-factor in $G$ such that $e_1,e_2$ are in the same cycle. If $F'$ contains $u_1u_4$ but not $u_2u_3$,
then $(F' \setminus \{u_1u_4\}) \cup \{ u_1v_1, v_1v_2, v_2v_3, v_3v_4, u_4v_4\}$ is a 2-factor in $G$
such that $e_1,e_2$ are contained in the same cycle. The same argument holds if $F'$ contains $u_2u_3$
but not $u_1u_4$. If $F'$ contains both, then $(F' \setminus \{u_1u_4, u_2u_3\}) \cup \{u_1v_1,u_2v_2,
u_3v_3,u_4v_4,v_1v_4,v_2v_3\}$ is a 2-factor in $G$ such that $e_1,e_2$ are in the same cycle of the
2-factor.
\hfill $\Box$

Another possible extension is to consider 2-factors including two specified edges, and characterize
those 2-edge-connected planar cubic graphs that have a pair of edges such that every 2-factor
containing both of them is a Hamiltonian cycle.

\end{document}